\documentclass[reqno,12pt]{amsart}
\usepackage[T1]{fontenc}
\usepackage{lmodern}
\usepackage{amsmath,amssymb,amsthm,mathrsfs}
\usepackage[margin=1in]{geometry}
\usepackage[expansion=false]{microtype}
\usepackage{graphicx,booktabs}
\usepackage{placeins}
\usepackage[
  colorlinks=true,
  linkcolor=blue,
  citecolor=green,
  urlcolor=blue
]{hyperref}

\newtheorem{theorem}{Theorem}[section]
\newtheorem{lemma}[theorem]{Lemma}
\newtheorem{proposition}[theorem]{Proposition}
\newtheorem{corollary}[theorem]{Corollary}
\theoremstyle{definition}

\theoremstyle{remark}
\newtheorem{remark}[theorem]{Remark}
\numberwithin{equation}{section}
\newcommand{\R}{\mathbb R}
\newcommand{\C}{\mathbb C}

\newcommand{\supp}{\operatorname{supp}}

\newcommand{\norm}[1]{\left\lVert#1\right\rVert}

\newcommand{\ip}[2]{\left\langle#1,#2\right\rangle}

\allowdisplaybreaks[4]

\makeatletter
\renewcommand{\@captionfont}{\normalfont\small}
\makeatother
\makeatletter
\def\@settitle{\begin{center}%
  \normalfont\Large\bfseries
  \@title
  \end{center}%
}
\makeatother
\title[Inverse random source problem]{Stochastic inverse source problem in inhomogeneous media: far-field statistics}
\author[D. Qiu]{Dong Qiu}
\address{School of Mathematical Sciences, Zhejiang University, Hangzhou 310058, China}
\email{qiudong@zju.edu.cn}
\author[X. Xu]{Xiang Xu}
\address{School of Mathematical Sciences, and Center for Interdisciplinary Applied Mathematics, Zhejiang University, Hangzhou 310058, China}
\email{xxu@zju.edu.cn}
\author[Y. Zhao]{Yue Zhao}
\address{School of Mathematics and Statistics, and Key Lab NAA-MOE, Central China Normal University, Wuhan 430079, China}
\email{zhaoy@ccnu.edu.cn}
\date{}

\newcommand{\E}{\mathbb E}
\newcommand{\HH}{\mathcal H}
\newcommand{\FF}{\mathcal F}
\newcommand{\im}{\operatorname{Im}}
\newcommand{\re}{\operatorname{Re}}
\newcommand{\sphere}{\mathbb S^2}
\renewcommand{\ip}[2]{\left\langle#1,#2\right\rangle_{\HH}}
\newcommand{\meanback}{\mathcal B}
\newcommand{\covback}{\mathcal R}
\newcommand{\ZZ}{\mathcal Z}
\newcommand{\YY}{\mathcal Y}

\begin{document}
\begin{abstract}
This paper is concerned with an inverse random source problem for the stochastic Helmholtz equation in a known three-dimensional inhomogeneous medium.
The source is driven by spatial white noise.
We seek to recover its mean and strength from the mean and cross-frequency covariance of the far-field pattern.
For a smooth nontrapping wave speed, we use the spectral measure of the acoustic operator to relate the far-field data to the source statistics.
Under a priori regularity bounds, this relation gives conditional stability for the source mean and the source strength.
We derive explicit reconstruction formulas for both quantities and bound the errors due to frequency truncation and data perturbations.
We also obtain mean-square error bounds for reconstructions from finitely many independent experiments.
Numerical experiments illustrate the reconstruction accuracy and the effects of bandwidth and sample size.
\end{abstract}
\maketitle

\section{Introduction}\label{sec:introduction}
Inverse source scattering problems arise in imaging and wave propagation, where the aim is to determine an unknown source from measurements of its radiated field.
When the sources exhibit random fluctuations, statistical information extracted from the measured fields can be used to infer source locations, spatial distributions, and statistical properties.
For example, optical fluctuation imaging uses temporal emission statistics to resolve fluorescent emitters~\cite{Dertinger2009}.
Microphone-array methods locate acoustic emissions by using source coherence and cross-spectral data~\cite{Sijtsma2007}.
Ambient seismic correlations provide information about distributed noise-source power in a heterogeneous Earth~\cite{Ermert2017}.
In this paper, we use far-field statistics to recover the mean and strength of a source driven by spatial white noise in a known inhomogeneous acoustic medium.

\subsection{Problem formulation}
Let $c\in C^\infty(\R^3)$ be a known wave speed satisfying
\begin{equation}\label{eq:medium}
0<c_-\le c(x)\le c_+<\infty,\qquad c(x)=1\quad (|x|\ge R_0).
\end{equation}
Assume that the Hamiltonian flow generated by $p(x,\xi)=c(x)^2|\xi|^2$ is nontrapping on the energy surface $p^{-1}(1)$, in the sense that every trajectory escapes every compact set as $t\to\pm\infty$.
Let $D\subset B_{R_0}$ be a fixed compact source region.
The random source is given by $f+\sqrt\mu\,\dot W$, with
\begin{equation}\label{eq:source}
\supp f\cup\supp\mu\subset D,\qquad \mu\ge0,
\end{equation}
where $W$ is real spatial white noise, $f\in L^2(D)$ is the mean of the random source, and $\mu\in L^\infty(D)$ is its strength.
The white noise is defined on a probability space $(\Omega,\mathscr A,\mathbb P)$, and $\E$ denotes expectation with respect to $\mathbb P$.
The notation $\dot W$ denotes the generalized random field associated with $W$, with $\langle\dot W,\varphi\rangle=W(\varphi)$ for real test functions $\varphi$.
The random component $\sqrt\mu\,\dot W$ has covariance kernel $\mu(x)\delta(x-y)$.

For~$k>0$, the radiated field satisfies
\begin{equation}\label{eq:pde}
\Delta u+k^2c^{-2}u=f+\sqrt\mu\,\dot W\quad\hbox{in }\R^3,
\qquad \lim_{r\to\infty}r(\partial_r u-iku)=0.
\end{equation}
We consider the outgoing mild solution of~\eqref{eq:pde} and denote its far-field pattern by $U(\theta,k)$, $\theta\in\sphere$.
Here $\sphere$ is the unit sphere, and all sphere integrals use the standard surface measure, so that $\int_{\sphere}1\,d\theta=4\pi$.
The far-field statistics are the mean $\E U(\theta,k)$ and the centered cross-frequency covariance
\begin{equation}\label{eq:statistics}
C(\theta,\eta;k,\ell)=\E\big((U(\theta,k)-\E U(\theta,k))
\overline{(U(\eta,\ell)-\E U(\eta,\ell))}\big).
\end{equation}
The same spatial noise realization is used at all frequencies within each experiment, and different experiments are independent.
The covariance $C$ is obtained by pairing the measurements at different frequencies.

The inverse problem is to reconstruct the source mean $f$ and strength $\mu$ from the far-field mean $\E U(\theta,k)$ and centered cross-frequency covariance $C(\theta,\eta;k,\ell)$ for all $\theta,\eta\in\sphere$ and all $k,\ell$ in a prescribed frequency interval.

\subsection{Previous work}    
Bao, Chen, and Li~\cite{BCL} established the existence and uniqueness of a mild solution to the stochastic Helmholtz equation.
They derived statistical integral equations and developed a regularized method to reconstruct the source statistics from multifrequency near-field data.
Frequency-averaged measurements for rough Gaussian sources were studied in~\cite{LHL} to obtain uniqueness results of the inverse random source problems for acoustic and elastic waves.
Stability for generalized Gaussian random sources in the presence of a potential was studied for the Helmholtz and Navier equations in~\cite{WXZ}.
For rough Gaussian sources, the recovery concerns the principal symbol of the covariance operator. In \cite{LLW2022FF}, the principal symbols of the covariance and relation operators were recovered from high-frequency averages of a single far-field realization in a homogeneous medium.

Multifrequency stability for deterministic sources was studied in~\cite{BLT,LZZ}.
In particular, \cite{LZZ} used spectral information, scattering resolvent estimates, and analytic continuation to establish stability in a known inhomogeneous medium.
Stability estimates for recovering white-noise source strengths from near-field correlations were established in~\cite{LL}.
Inverse random potential problems, in which the randomness enters a multiplicative coefficient, were studied for two-dimensional Schr\"odinger scattering~\cite{LPS2008} and single-realization backscattering~\cite{CHL2019}.

In this work, we study the direct and inverse random source scattering problems in a known inhomogeneous medium.
We establish the existence and uniqueness of the outgoing mild solution and derive integral equations relating its far-field statistics to the source mean and strength.
The radiation-flux identity and Stone's formula give a spectral representation of the acoustic operator in terms of the far-field pattern.
Under this representation, the difference of two cross-frequency covariances determines the compression of a multiplication operator to the measured frequency band.
For a smooth nontrapping wave speed, we establish conditional stability in $L^2$ for the source mean and in $L^\infty$ for the source strength under a priori regularity bounds.
The estimates imply uniqueness from a fixed frequency interval and quantify the dependence of the regularity remainders on the bandwidth.

We further derive explicit reconstruction formulas for the mean and strength from the same spectral representation.
For the source mean, we derive a weighted backprojection formula from the adjoint far-field transform.
For the source strength, we fix a reference frequency and express the cross-frequency covariance as the far-field pattern of the product of the strength and a known reference kernel.
Applying the spectral identity and integrating over the reference directions give a positive weight and a weighted backprojection formula for the strength.
The reconstruction uses covariance data with one frequency fixed and the other varying over the measured band.
The resulting reconstructions converge in $L^2(D)$ as the bandwidth tends to infinity.
Under a graph-norm bound on $c^2f$ and a Sobolev bound on $\mu$, we estimate the reconstruction errors due to frequency truncation and data perturbations.

To reconstruct the source functions from repeated experiments, we replace the exact far-field statistics by the sample mean and cross-frequency sample covariance.
The observation and reference measurements in each pair share the same realization of the spatial white noise.
We derive finite-sample error estimates for both reconstructions in terms of the bandwidth and the number of experiments.
The numerical algorithm evaluates the truncated weighted backprojections by quadrature.
Numerical experiments with relative measurement noise demonstrate the reconstruction of smooth sources and illustrate the effects of bandwidth and sample size.

The outline of the paper is as follows.
Section~\ref{sec:direct} presents the outgoing mild solution and the spectral representation of the far-field statistics.
Section~\ref{sec:stability} is devoted to the stability estimates for the source mean and strength.
The reconstruction formulas and finite-sample error estimates are derived in Section~\ref{sec:reconstruction}.
Numerical experiments are presented in Section~\ref{sec:numerics}.
The paper is concluded in Section~\ref{sec:conclusion}.

\section{Forward problem and spectral representation}\label{sec:direct}
\subsection{Mild solution and far-field statistics}
We discuss the outgoing mild solution and derive its far-field statistics.
For real $\varphi,\psi\in L^2(D)$, white noise is a centered Gaussian linear map satisfying
\begin{equation}\label{eq:white}
\E(W(\varphi)W(\psi))=\int_D\varphi(y)\psi(y)\,dy.
\end{equation}
Extending this map complex linearly, we obtain for complex integrands
\begin{equation}\label{eq:isometry}
\E\left(\int_D v\,dW\ \overline{\int_D w\,dW}\right)
=\int_Dv\overline w\,dy.
\end{equation}

Let $G(x,y;k)$ be the outgoing kernel satisfying
\begin{equation}\label{eq:green}
(\Delta_x+k^2c(x)^{-2})G(x,y;k)=\delta_y.
\end{equation}
With this sign convention, $G(x,y;k)=-e^{ik|x-y|}/(4\pi|x-y|)$ when $c=1$.
The outgoing kernel has the far-field expansion
\begin{equation}\label{eq:kernelff}
G(r\theta,y;k)=r^{-1}e^{ikr}\Gamma(\theta,y;k)+O(r^{-2}),\qquad \theta\in\sphere.
\end{equation}
Following~\cite{BCL}, the mild solution is defined by
\begin{equation}\label{eq:mild}
u(x,k)=\int_DG(x,y;k)f(y)\,dy+
\int_DG(x,y;k)\sqrt{\mu(y)}\,dW_y.
\end{equation}
The stochastic integral is understood in $L^2(\Omega;L^2_{\mathrm{loc}}(\R^3))$.
We denote the stochastic integral in~\eqref{eq:mild} by $u^{\rm s}$.

\begin{proposition}\label{prop:direct}
For each $k>0$, the stochastic Helmholtz equation has a unique outgoing mild solution in $L^2(\Omega;L^2_{\mathrm{loc}}(\R^3))$, given by~\eqref{eq:mild}.
The outgoing condition holds in mean square on exterior spheres, and the solution satisfies
\begin{equation}\label{eq:msfarfield}
\left\|u(r\cdot,k)-\frac{e^{ikr}}rU(\cdot,k)\right\|_{L^2(\Omega;L^2(\sphere))}
=O(r^{-2}),
\end{equation}
where the far-field pattern $U(\theta,k)$ is given by
\begin{equation}\label{eq:randomff}
U(\theta,k)=\int_D\Gamma(\theta,y;k)f(y)\,dy+
\int_D\Gamma(\theta,y;k)\sqrt{\mu(y)}\,dW_y.
\end{equation}
The first moment and centered covariance are given by
\begin{align}
 \E U(\theta,k)&=\int_D\Gamma(\theta,y;k)f(y)\,dy,\label{eq:mean}\\
 C(\theta,\varphi;k,\ell)&=\E\big((U(\theta,k)-\E U(\theta,k))
 \overline{(U(\varphi,\ell)-\E U(\varphi,\ell))}\big)\notag\\
 &=\int_D\mu(y)\Gamma(\theta,y;k)
 \overline{\Gamma(\varphi,y;\ell)}\,dy.\label{eq:cov}
\end{align}
\end{proposition}

\begin{proof}
Fredholm theory and uniqueness for the homogeneous outgoing problem give the outgoing Green operator for each $k>0$.
For a bounded set $B$, define
\[
(K_Bv)(x)=\int_DG(x,y;k)\sqrt{\mu(y)}v(y)\,dy.
\]
In three dimensions, the Green function has a singularity of order $|x-y|^{-1}$, so
\begin{equation}\label{eq:hilbertschmidt}
\|K_B\|_{\rm HS}^2=\int_B\int_D|G(x,y;k)|^2\mu(y)\,dy\,dx<\infty.
\end{equation}
Here $\|\cdot\|_{\rm HS}$ denotes the Hilbert--Schmidt norm, applied to $K_B:L^2(D)\to L^2(B)$.
Let $(e_j)$ be a real orthonormal basis of $L^2(D)$.
Since $W(e_j)$ are independent Gaussian variables with variance one, we have
\begin{equation}\label{eq:mildseries}
u^{\rm s}|_B=\sum_{j=1}^\infty W(e_j)K_Be_j
\quad\hbox{in }L^2(\Omega;L^2(B)),\qquad
\E\|u^{\rm s}\|_{L^2(B)}^2=\|K_B\|_{\rm HS}^2.
\end{equation}
By the white-noise isometry, the construction is independent of the basis and consistent on nested bounded sets.
Adding the deterministic Green integral yields~\eqref{eq:mild}.
Testing the kernel equation and using stochastic Fubini give
\[
\int_{\R^3}u(\Delta+k^2c^{-2})\varphi\,dx
=\int_D f\varphi\,dx+W(\sqrt\mu\,\varphi),\qquad
\varphi\in C_c^\infty(\R^3).
\]
Thus~\eqref{eq:pde} holds in the sense of random distributions.

The asymptotic expansion of the outgoing kernel is uniform for $y\in D$.
The white-noise isometry gives
\[
\E\left\|u(r\cdot,k)-\frac{e^{ikr}}rU(\cdot,k)\right\|_{L^2(\sphere)}^2
\le Cr^{-4}(\|f\|_{L^1(D)}^2+\|\mu\|_{L^1(D)}).
\]
The far-field pattern is therefore given by~\eqref{eq:randomff}.
Taking expectations and applying the isometry give the first-moment and covariance identities.

The Cauchy data of $G(\cdot,y;k)$ on a fixed sphere enclosing $B_{R_0}$ are uniformly bounded for $y\in D$.
Applying $\partial_r-ik$ to the free-space kernel and its boundary normal derivative in the exterior Green representation gives
\[
(\partial_r-ik)G(r\theta,y;k)=O(r^{-2})
\]
uniformly for $\theta\in\sphere$ and $y\in D$.
At each exterior radius, the radial difference quotients of $G$ converge uniformly on $\sphere\times D$, so the white-noise isometry justifies differentiating~\eqref{eq:mild} in mean square.
The kernel estimate and the white-noise isometry yield
\[
\E\big\|r(\partial_r-ik)u(r\cdot,k)\big\|_{L^2(\sphere)}^2
\le Cr^{-2}\big(\|f\|_{L^1(D)}^2+\|\mu\|_{L^1(D)}\big),
\]
which proves the mean-square radiation condition.

For uniqueness, pair the difference of two outgoing solutions with any scalar random variable in $L^2(\Omega)$.
The pairing is a deterministic homogeneous outgoing solution, which vanishes by Rellich's lemma and unique continuation.
These pairings separate $L^2(\Omega;L^2(B))$, so the two solutions coincide. The proof is complete.
\end{proof}

Equations~\eqref{eq:mean} and~\eqref{eq:cov} are linear integral equations for the source statistics, as in~\cite{BCL}.
The equations for $f$ and $\mu$ are decoupled because centering removes the unknown mean from the covariance equation.
We analyze the two equations using separate regularity assumptions and data norms.

The Green representation and reciprocity follow from Green's identities for inhomogeneous acoustic scattering, see~\cite[Chapter~8]{CK2019}.
The Green kernel satisfies bilinear reciprocity, $G(x,y;k)=G(y,x;k)$.
Under the normalization~\eqref{eq:green}, mixed reciprocity gives
\begin{equation}\label{eq:mixed}
\Gamma(\theta,y;k)=-\frac1{4\pi}v(y,-\theta;k),
\end{equation}
where $v$ denotes the total field corresponding to the incident plane wave $e^{-ik\theta\cdot y}$ in the same background medium.

\subsection{Radiation flux and the spectral measure}
Set
\begin{equation}\label{eq:operator}
\HH=L^2(\R^3,c^{-2}dx),
\end{equation}
and realize $-c^2\Delta$ as the nonnegative self-adjoint operator associated with the closed form $\int\nabla v\cdot\overline{\nabla w}\,dx$ on $H^1(\R^3)$.
Throughout the paper, inner products are linear in the first argument.

Let $R_+(k)=(-c^2\Delta-(k+i0)^2)^{-1}$ denote the outgoing resolvent acting on compactly supported data.
The notation $k+i0$ denotes the boundary value obtained from $k+i\varepsilon$ as $\varepsilon\downarrow0$.
\begin{lemma}\label{lem:spectrum}
The spectrum of $-c^2\Delta$ is purely absolutely continuous and equals $\{\lambda\in\R:\lambda\ge0\}$.
The operator $-c^2\Delta$ has no eigenvalues.
\end{lemma}
\begin{proof}
Nonnegativity excludes negative eigenvalues.
If $-c^2\Delta v=0$, the form identity gives $\int|\nabla v|^2dx=0$, and $v\in\HH$ implies $v=0$.
An $L^2$ eigenfunction associated with a positive eigenvalue satisfies the free Helmholtz equation outside a ball.
Exterior $L^2$ uniqueness and unique continuation then force it to vanish.
On each compact positive energy interval, Proposition~\ref{prop:cutoffresolvent} gives norm-continuous boundary values of the cutoff outgoing resolvent.
Stone's formula then gives a locally continuous density for the spectral measure of each compactly supported smooth vector.
By density, there is no singular continuous spectrum on $(0,\infty)$.
The trivial kernel of $-c^2\Delta$ excludes a spectral atom at zero.
Finally, exterior Weyl sequences show that every nonnegative energy belongs to the spectrum.
\end{proof}

For deterministic forcing $g\in\HH$, Green's identity and Stone's formula relate the spectral measure of $-c^2\Delta$ to the far-field energy.
Integration over all positive frequencies gives the full source energy.
Integration over a finite frequency band gives the energy of the source projected onto that band.
For the generalized plane-wave representation and the spectral theory of self-adjoint scattering operators, see~\cite[Theorem~4.20]{DZ} and~\cite{RS}.
We derive the identity under the far-field convention~\eqref{eq:kernelff}.
Write $E_{-c^2\Delta}(\lambda)=\mathbf1_{[0,\lambda]}(-c^2\Delta)$, $\lambda\ge0$, for the spectral family of $-c^2\Delta$.

\begin{proposition}\label{prop:flux}
For compactly supported $g\in\HH$, let $u_g=-R_+(k)g$ with far field $u_g^\infty$.
Then
\begin{align}
d\ip{E_{-c^2\Delta}(k^2)g}{g}
&=\frac2\pi k^2\|u_g^\infty(\cdot,k)\|_{L^2(\sphere)}^2\,dk,
\label{eq:spectralmeasure}\\
\|g\|_{\HH}^2
&=\frac2\pi\int_0^\infty k^2\|u_g^\infty(\cdot,k)\|_{L^2(\sphere)}^2\,dk.
\label{eq:flux}
\end{align}
Here the differential denotes the spectral measure in the frequency variable $k>0$.
\end{proposition}
\begin{proof}
Set $v=R_+(k)g=-u_g$.
Integrating by parts in $B_\rho$, we obtain
\[
\int_{B_\rho}(|\nabla v|^2-k^2c^{-2}|v|^2)dx-
\int_{\partial B_\rho}\partial_r v\,\overline v\,ds
=\int_{B_\rho}c^{-2}g\overline v\,dx.
\]
Taking imaginary parts and using the inner-product convention, we obtain
\[
\im\ip{R_+(k)g}{g}
=\lim_{\rho\to\infty}\im\int_{\partial B_\rho}\partial_r v\,\overline v\,ds
=k\norm{u_g^\infty}_{L^2(\sphere)}^2.
\]
Applying Stone's formula~\cite[Theorem~B.10]{DZ} with $\lambda=k^2$ and $d\lambda=2k\,dk$ yields~\eqref{eq:spectralmeasure}.
Integrating over the spectrum gives~\eqref{eq:flux} by Lemma~\ref{lem:spectrum}, with monotone convergence at zero and infinity.
\end{proof}

For compactly supported data, define the far-field transform by
\begin{equation}\label{eq:transform}
(\FF g)(k\theta)=-\sqrt{2/\pi}\,u_g^\infty(\theta,k).
\end{equation}
For $L>0$, set $\Pi_L=\mathbf1_{\{0\le\lambda\le L^2\}}(-c^2\Delta)$ and $\chi_L=\mathbf1_{B_L}$.
In operator expressions, $I$ denotes the identity operator on the relevant Hilbert space.
\begin{corollary}\label{cor:transform}
The map $\FF$ extends to an isometry $\HH\to L^2(\R^3)$, and
\begin{equation}\label{eq:intertwine}
\FF^*\FF=I,\qquad \FF^*\chi_L\FF=\Pi_L,\qquad
\chi_L\FF=\FF\Pi_L.
\end{equation}
\end{corollary}
\begin{proof}
The full-frequency identity and density give the isometry.
Integrating the spectral identity up to $L$ and polarizing, we obtain $\FF^*\chi_L\FF=\Pi_L$.
Taking squared norms then gives $\chi_L\FF(I-\Pi_L)=0$ and $(I-\chi_L)\FF\Pi_L=0$.
The last assertion follows by combining these identities.
\end{proof}

\subsection{Covariance as a compressed multiplication operator}
Consider two random sources satisfying the assumptions of the problem formulation.
Let $\E U_j,C_j$ denote the far-field statistics of source $j$ and write $\mathcal C=C_1-C_2$.
For $K>0$, identify $\xi=k\theta$, $\zeta=\ell\varphi$, and use Lebesgue measure $d\xi=k^2dk\,d\theta$ in frequency space.
Define the covariance data operator by
\begin{equation}\label{eq:datamu}
(A_Kb)(\xi)=\frac2\pi\int_{B_K}\mathcal C(\xi,\zeta)b(\zeta)\,d\zeta,
\qquad \varepsilon_K=\norm{A_K}_{L^2(B_K)\to L^2(B_K)}.
\end{equation}

In $\HH$, the covariance form of the weighted random source $c^2\sqrt\mu\,\dot W$ is represented by multiplication by $c^2\mu$.
Applying the isometric far-field transform and restricting both frequencies to the measured band give the following compression identity.
Set $a=c^2(\mu_1-\mu_2)$ and let $M_a$ denote multiplication by $a$ on $\HH$.
\begin{lemma}\label{lem:compression}
Extend $A_K$ by zero to $L^2(\R^3)$.
Then
\begin{equation}\label{eq:compression}
 A_K=\FF\Pi_KM_{c^2(\mu_1-\mu_2)}\Pi_K\FF^*,\qquad
 \varepsilon_K=\norm{\Pi_KM_{c^2(\mu_1-\mu_2)}\Pi_K}_{\HH\to\HH}.
\end{equation}
Moreover,
\begin{equation}\label{eq:schur}
\varepsilon_K\le\frac83 K^3
\sup_{\theta,\varphi\in\sphere,\ 0<k,\ell\le K}|\mathcal C(\theta,\varphi;k,\ell)|.
\end{equation}
\end{lemma}
\begin{proof}
The weighted random source $c^2\sqrt{\mu_j}\dot W$ is a generalized random field on $\HH$.
Its pairing with $v\in\HH$ is $\int\sqrt{\mu_j}\,\overline v\,dW$, and the covariance form on $\HH$ is represented by multiplication by $c^2\mu_j$.
Let $X_j$ denote the stochastic spectral coefficient $-\sqrt{2/\pi}(U_j-\E U_j)$, understood as a distribution.
For spectral test functions $b,v$ on compact positive frequency annuli, the white-noise isometry and~\eqref{eq:transform} give
\[
\E\left(\left(\int X_j(\xi)\overline{b(\xi)}d\xi\right)
\overline{\left(\int X_j(\zeta)\overline{v(\zeta)}d\zeta\right)}\right)
=\ip{M_{c^2\mu_j}\FF^*v}{\FF^*b}.
\]
Subtracting the two covariance identities, we obtain $A_K=\chi_K\FF M_a\FF^*\chi_K$ in the weak operator sense.
Since $\mu_1-\mu_2\in L^\infty(D)$, the multiplier is bounded, and density extends the identity to all $L^2$ test functions.
Applying~\eqref{eq:intertwine} gives the operator identity.
For every bounded operator $T$ on $\HH$, the isometry gives $\norm{\FF T\FF^*}=\norm T$, proving the norm equality.
Schur's test gives $\varepsilon_K\le(2/\pi)|B_K|\sup|\mathcal C|$, which proves~\eqref{eq:schur}.
\end{proof}

\begin{corollary}\label{cor:fullbandnorms}
The far-field statistics satisfy the norm identities
\begin{align}
\|c(f_1-f_2)\|_{L^2}^2
&=\frac2\pi\int_0^\infty k^2\|\E U_1(\cdot,k)-\E U_2(\cdot,k)\|_{L^2(\sphere)}^2\,dk,
\label{eq:fullbandmean}\\
\|c^2(\mu_1-\mu_2)\|_{L^\infty}
&=\sup_{K>0}\varepsilon_K.\label{eq:fullbandnorms}
\end{align}
\end{corollary}
\begin{proof}
The first identity follows directly from~\eqref{eq:flux}.
For the second identity, $\Pi_KM_a\Pi_K\to M_a$ strongly as $K\to\infty$, and each compressed operator has norm at most $\|M_a\|=\|a\|_{L^\infty}$.
The reverse inequality follows by applying strong convergence to each unit vector and taking the supremum.
Together with~\eqref{eq:compression}, this proves the identity.
\end{proof}

\section{Inverse source stability}\label{sec:stability}
In this section, we establish stability estimates for the source mean and strength.
Let $r>0$, $s>3/2$, and fix the a priori bounds $Q,M>0$.
The fractional powers of $I-c^2\Delta$ are defined by the spectral calculus of $-c^2\Delta$.
We equip $D((I-c^2\Delta)^{r/2})$ with the equivalent graph norm $\|(I-c^2\Delta)^{r/2}v\|_{\HH}$, where $D(T)$ denotes the domain of an operator $T$.
For $j=1,2$, assume that
\begin{align}
 &\supp f_j\subset D,\qquad
 \norm{(I-c^2\Delta)^{r/2}(c^2 f_j)}_{\HH}\le Q,\label{eq:priorf}\\
 &\mu_j\ge0,\quad\supp\mu_j\subset D,\qquad
 \norm{\mu_j}_{H^s(\R^3)}\le M.\label{eq:priormu}
\end{align}
The source functions are extended by zero to the whole space.
Choose
\begin{equation}\label{eq:gamma}
0<\gamma<\min(1,s-3/2).
\end{equation}

Fix $0<k_0<k_1$ and let $I=\{k:k_0\le k\le k_1\}$ be the observation interval.
Suppose that the far-field data differences on $I$ satisfy
\begin{equation}\label{eq:small}
\sup_{k\in I} k\norm{\E U_1(\cdot,k)-\E U_2(\cdot,k)}_{L^2(\sphere)}\le Q\delta_f,
\qquad
\sup_{\theta,\varphi\in\sphere,\ k,\ell\in I}|\mathcal C(\theta,\varphi;k,\ell)|\le M\delta_\mu.
\end{equation}
Here $\delta_f$ and $\delta_\mu$ are normalized bounds for the data differences.

Define the weighted norm $D_K$ of the difference between the far-field means by
\begin{equation}\label{eq:dataf}
D_K^2=\int_0^K k^2\norm{\E U_1(\cdot,k)-\E U_2(\cdot,k)}_{L^2(\sphere)}^2\,dk.
\end{equation}

In the following theorem, we present conditional stability estimates for the stochastic inverse source problem using the multi-frequency far-field patterns.
\begin{theorem}\label{thm:fixed}
Assume the standing hypotheses on $c,D$ and~\eqref{eq:priorf}--\eqref{eq:small}.
There exist constants $C>0$ and $\delta_0\in(0,e^{-e})$, independent of the sources, $Q,M$, and the data differences, such that
\begin{align}
 \norm{f_1-f_2}_{L^2(\R^3)}^2&\le C Q^2 (\ln|\ln\delta_f|)^{-2r},\label{eq:fixedf}\\
 \norm{\mu_1-\mu_2}_{L^\infty(\R^3)}&\le C M (\ln|\ln\delta_\mu|)^{-\gamma},\label{eq:fixedmu}
\end{align}
with the first estimate valid for $0<\delta_f\le\delta_0$ and the second for $0<\delta_\mu\le\delta_0$.
The constants may depend on $c,D,I,r,s$, and $\gamma$.
In particular, if the far-field means agree on $\sphere\times I$ and the centered cross-frequency covariances agree on $\sphere\times\sphere\times I\times I$, then $f_1=f_2$ and $\mu_1=\mu_2$.
\end{theorem}

Moreover, we prove increasing stability estimates for the stochastic inverse source problem.
Each estimate consists of a data discrepancy term and a remainder controlled by the a priori regularity bounds.
The remainders tend to zero as the bandwidth increases.
\begin{theorem}\label{thm:increasing}
Let $0<\delta_f,\delta_\mu\le\delta_0$ and assume the hypotheses of Theorem~\ref{thm:fixed}.
Suppose also that the mean and cross-frequency covariance are available for $0<k,\ell<K$, where $K\ge\max(1,k_1)$.
Then
\begin{align}
 \norm{f_1-f_2}_{L^2}^2
 &\le C\left(D_K^2+\frac{Q^2}{\max(K,\ln|\ln\delta_f|)^{2r}}\right),\label{eq:incf}\\
 \norm{\mu_1-\mu_2}_{L^\infty}
 &\le C\left(\varepsilon_K+\frac{M}{\max(K,\ln|\ln\delta_\mu|)^{\gamma}}\right).\label{eq:incmu}
\end{align}
The constant is also independent of $K$.
\end{theorem}

\subsection{Stability of the source mean}
Let $q_f=f_1-f_2$ be the difference of the source means and set $g=c^2q_f$.
Then $u_g^\infty=\E U_1-\E U_2$ and $\norm{(I-c^2\Delta)^{r/2}g}_{\HH}\le2Q$.
Functional calculus and~\eqref{eq:flux} give
\begin{equation}\label{eq:meantail}
\int_L^\infty k^2\norm{\E U_1(\cdot,k)-\E U_2(\cdot,k)}_{L^2(\sphere)}^2dk
\le C Q^2(1+L^2)^{-r}.
\end{equation}
Using $\norm{g}_{\HH}^2=\int c^2|q_f|^2dx\asymp\norm{q_f}_{L^2}^2$ and splitting~\eqref{eq:flux} at $L\ge1$, we obtain
\begin{equation}\label{eq:finitef}
\norm{q_f}_{L^2}^2\le C(D_L^2+Q^2L^{-2r}).
\end{equation}

The following resolvent estimates allow analytic continuation beyond the observation interval.

\begin{proposition}\label{prop:cutoffresolvent}
Assume~\eqref{eq:medium} and the nontrapping condition.
The outgoing resolvent $R(\lambda)=(-c^2\Delta-\lambda^2)^{-1}$ admits a meromorphic continuation from the physical upper half-plane to $\C$ as a map $\HH_{\rm comp}\to H^2_{\rm loc}(\R^3)$.
Here $\lambda$ is a complex frequency, $\HH_{\rm comp}$ consists of the compactly supported elements of $\HH$, and $R(k+i0)=R_+(k)$ for $k>0$.
There exist $\eta,d_0>0$ such that this continuation has no poles in
\begin{equation}\label{eq:strip}
\re\lambda>-\eta,\qquad |\im\lambda|<d_0.
\end{equation}
For each fixed $\chi\in C_c^\infty(\R^3)$, there exists $C_\chi>0$ such that, throughout the strip~\eqref{eq:strip},
\begin{equation}\label{eq:cutoffresolvent}
\|\chi R(\lambda)\chi\|_{\HH\to\HH}
\le\frac{C_\chi}{1+|\lambda|},\qquad
\|\chi R(\lambda)\chi\|_{\HH\to H^2}
\le C_\chi(1+|\lambda|).
\end{equation}
\end{proposition}
\begin{proof}
The self-adjoint form realization of $-c^2\Delta$ has domain $H^2(\R^3)$.
Outside $B_{R_0}$, it agrees with $-\Delta$, and elliptic regularity and the Rellich embedding give compactness of the interior part of $(-c^2\Delta+i)^{-1}$.
Thus $-c^2\Delta$ satisfies Definition~4.1 of~\cite{DZ}.
Theorem~4.4 of that reference, in particular equation~(4.2.5), gives meromorphic continuation to $\C$ in dimension three.

Equip $\R^3$ with the metric $\mathfrak g_{ij}=c^{-2}\delta_{ij}$ and volume $dV_{\mathfrak g}=c^{-3}dx$.
We write $\Delta_{\mathfrak g}$ for the Laplace--Beltrami operator of this metric.
The unitary map $\mathscr U:\HH\to L^2(dV_{\mathfrak g})$, defined by $\mathscr Uv=c^{1/2}v$, satisfies
\begin{equation}\label{eq:metricconjugation}
\mathscr U(-c^2\Delta)\mathscr U^{-1}=-\Delta_{\mathfrak g}+V_c,
\qquad V_c=\tfrac12c\Delta c-\tfrac34|\nabla c|^2.
\end{equation}
The metric is a smooth, compactly supported perturbation of the Euclidean metric, and $V_c\in C_c^\infty(\R^3;\R)$.
The geodesic Hamiltonian is $p(x,\xi)=c^2|\xi|^2$.
Homogeneity in $\xi$ extends nontrapping at energy one to every compact positive energy interval.

Let $\hbar>0$ be a semiclassical parameter.
Denote by $R_{0,\hbar}(z)=(-\hbar^2\Delta_{\mathfrak g}-z)^{-1}$ the resolvent initially defined for $\im z>0$ and continued with the outgoing condition.
We apply~\cite[Theorem~6.10]{DZ} on $\{E:1/2\le E\le2\}$.
For any fixed $a>0$, the outgoing resolvent $R_{0,\hbar}(z)$ has no poles for $1/2\le\re z\le2$, $|\im z|\le a\hbar$, and satisfies
\begin{equation}\label{eq:semiclassicalinput}
\|\chi R_{0,\hbar}(z)\chi\|_{L^2(dV_{\mathfrak g})\to L^2(dV_{\mathfrak g})}
\le C\hbar^{-1}
\end{equation}
for sufficiently small $\hbar$.
The semiclassical Sobolev space $H^2_\hbar$ is $H^2(\R^3)$ equipped with the norm
\[
\|v\|_{H^2_\hbar}=\|(I-\hbar^2\Delta_{\mathfrak g})v\|_{L^2(dV_{\mathfrak g})}.
\]
For the same cutoff resolvent, the theorem gives an $L^2(dV_{\mathfrak g})\to H^2_\hbar$ operator norm bound of $C\hbar^{-1}$.

The additional potential is $\hbar^2V_c$.
Choose $\chi=1$ near its support.
On a ball containing the support, $\|\hbar^2\chi R_{0,\hbar}(z)V_c\|\le C\hbar$, so $I+\hbar^2\chi R_{0,\hbar}(z)V_c$ is invertible by a Neumann series.
The resolvent identity gives the same pole-free region and $C\hbar^{-1}$ cutoff $L^2$ estimate for $\hbar^2\mathscr U(-c^2\Delta)\mathscr U^{-1}$, including the lower-order term in~\eqref{eq:metricconjugation}.

Write $\lambda=t+i\sigma$, choose $\hbar=t^{-1}$, and set $z=\hbar^2\lambda^2=1-\hbar^2\sigma^2+2i\hbar\sigma$.
For bounded $|\sigma|$ and large $t$, this parameter lies in the preceding semiclassical region.
Rescaling the resolvent introduces the factor $\hbar^2$, and hence $\|\chi R(\lambda)\chi\|_{\HH\to\HH}\le C/t$.
Applying interior elliptic estimates with a larger cutoff $\chi_1\in C_c^\infty(\R^3)$ satisfying $\chi_1=1$ near $\supp\chi$, we obtain
\[
\|\chi R(\lambda)\chi f\|_{H^2}
\le C\big(\|f\|_{\HH}
+(1+|\lambda|^2)\|\chi_1R(\lambda)\chi f\|_{\HH}\big)
\le C(1+|\lambda|)\|f\|_{\HH}.
\]
This proves~\eqref{eq:cutoffresolvent} at high frequency in a fixed strip.

It remains to consider bounded frequencies.
At a positive real frequency, Green's identity gives zero radiated flux for a homogeneous outgoing solution, which vanishes by Rellich's lemma and unique continuation.
A Laurent coefficient at a real pole would give a nonzero homogeneous outgoing solution, so no positive real poles exist.
At zero, set $q=c^{-2}-1$ and $\Phi_\lambda(x,y)=e^{i\lambda|x-y|}/(4\pi|x-y|)$.
On a ball containing $\supp q$, the Lippmann--Schwinger operator is $I-\lambda^2\Phi_\lambda q$.
Here $\Phi_\lambda q$ denotes the integral operator with kernel $\Phi_\lambda(x,y)q(y)$ on that ball.
The integral part is holomorphic and vanishes at zero, so the operator is invertible and the local resolvent is holomorphic there.
Meromorphy and the absence of real poles give a pole-free neighborhood of the remaining compact real-frequency interval.
Combining this neighborhood with the high-frequency strip and decreasing the width if necessary, we obtain~\eqref{eq:strip}.
Continuity on the remaining compact set completes the proof of~\eqref{eq:cutoffresolvent}.
\end{proof}

\begin{lemma}\label{lem:strip}
The kernel $\Gamma(\theta,\cdot;\lambda)$ extends holomorphically to the strip~\eqref{eq:strip} with values in $C^0(D)$, and
\begin{equation}\label{eq:kernelbound}
\sup_{\theta\in\sphere,\ y\in D}|\Gamma(\theta,y;\lambda)|
\le C(1+|\lambda|)^3.
\end{equation}
\end{lemma}
\begin{proof}
Let $q=c^{-2}-1$ and choose cutoffs equal to one near $D\cup\supp q$.
Mixed reciprocity~\eqref{eq:mixed} gives
\begin{equation}\label{eq:totalfieldresolvent}
\Gamma(\theta,y;\lambda)=-\frac1{4\pi}
\big(e^{-i\lambda\theta\cdot y}
+\lambda^2R(\lambda)(c^2q e^{-i\lambda\theta\cdot(\cdot)})(y)\big).
\end{equation}
The identity extends from positive real $\lambda$ to the strip by analytic continuation.
The compactly supported source satisfies $\sup_\theta\|c^2q e^{-i\lambda\theta\cdot(\cdot)}\|_{\HH}\le C e^{d_0R_0}$ in the strip.
By the second estimate in~\eqref{eq:cutoffresolvent}, the outgoing correction is bounded in $H^2$ near $D$ by $C(1+|\lambda|)^3$.
The embedding $H^2\hookrightarrow C^0$ in dimension three gives~\eqref{eq:kernelbound}, uniformly in $\theta$ and $y$.
The cutoff resolvent is holomorphic with values in $H^2_{\rm loc}$, and the same embedding gives holomorphy with values in $C^0(D)$.
The proof is complete.
\end{proof}

\begin{lemma}\label{lem:propagation}
There exist $a_0,c_0>0$ and a function $\omega:\{k\ge0\}\to\{t:0<t\le1\}$ with $\omega(k)\ge c_0e^{-a_0k}$ such that every holomorphic $H$ in the strip~\eqref{eq:strip} satisfying $|H|\le1$ there and $|H|\le\delta<1$ on $I$ obeys
\begin{equation}\label{eq:propagation}
|H(k)|\le\delta^{\omega(k)}\quad(k\ge0).
\end{equation}
One may take $\omega=1$ on $I$.
\end{lemma}
\begin{proof}
In the strip~\eqref{eq:strip} with $I$ removed, consider the harmonic measure of the slit with boundary values one on the slit and zero elsewhere.
Construct it by exhaustion with bounded truncated slit domains and zero values on the artificial boundaries.
The harmonic measure is positive in the domain and tends to one at the regular slit endpoints.
Fix $k_c>k_1+d_0$.
For $k\ge k_c$, Harnack's inequality along a chain of overlapping fixed-radius disks gives $\omega(k)\ge C_H^{-J(k)}\omega(k_c)$, where $J(k)\le C_1+C_2k$.
Positivity and regularity at the slit endpoints give a uniform positive lower bound on the remaining compact real intervals.
These bounds prove the exponential lower bound.

Apply the two-constants principle to the bounded-above subharmonic function $\log|H|$.
On each bounded exhaustion domain, $\log|H|\le0$ on the artificial boundary, so comparison with the harmonic measure multiplied by $\log\delta$ applies.
Passing to the limit gives $\log|H(k)|\le\omega(k)\log\delta$, which proves the estimate.
\end{proof}

It follows from~\eqref{eq:mean} and Lemma~\ref{lem:strip} that $\E U_1(\cdot,z)-\E U_2(\cdot,z)$ is holomorphic with values in $L^2(\sphere)$ and satisfies $\norm{\E U_1(\cdot,z)-\E U_2(\cdot,z)}\le CQ(1+|z|)^3$.
We continue the weighted squared $L^2(\sphere)$ norm of the difference between the far-field means by setting
\begin{equation}\label{eq:energy}
F(z)=z^2\int_{\sphere}\big(\E U_1(\theta,z)-\E U_2(\theta,z)\big)\overline{\E U_1(\theta,\overline z)-\E U_2(\theta,\overline z)}\,d\theta.
\end{equation}
Since the strip is invariant under complex conjugation, the reflected conjugate in the second factor is holomorphic in $z$.
The bilinear integral therefore defines a holomorphic function $F$.
On the real axis, $F(k)=k^2\norm{\E U_1(\cdot,k)-\E U_2(\cdot,k)}^2\ge0$, also for complex-valued source means.

\begin{proposition}\label{prop:meancont}
There are $a,b,\beta,C>0$ such that
\begin{equation}\label{eq:meancont}
F(k)\le CQ^2\exp\big(bk-\beta|\ln\delta_f|e^{-ak}\big),\qquad k\ge0,
\end{equation}
provided~\eqref{eq:small} holds with $0<\delta_f<1$.
\end{proposition}
\begin{proof}
Fix $b>0$ and take $B\ge1$ sufficiently large.
The polynomial growth bound and the lower bound on $\re z$ show that $H(z)=e^{-bz}F(z)/(BQ^2)$ has modulus at most one throughout the strip.
On $I$, the data bound and $k>0$ give a modulus of at most $\delta_f^2$.
Applying Lemma~\ref{lem:propagation} and restoring the exponential factor prove the estimate with $a=a_0$ and $\beta=2c_0$.
\end{proof}

\begin{proof}[Proof of~\eqref{eq:fixedf} and~\eqref{eq:incf}]
Set $t=|\ln\delta_f|$ and choose
\begin{equation}\label{eq:choicef}
L_f=\frac1{2a}\ln t,
\end{equation}
where $a$ is the constant in Proposition~\ref{prop:meancont}.
Take $\delta_0$ sufficiently small that $L_f\ge1$.
Integrating~\eqref{eq:meancont} from zero to $L_f$, we obtain
\[
D_{L_f}^2\le CQ^2L_f\exp(bL_f-\beta t^{1/2})
=CQ^2L_f t^{b/(2a)}e^{-\beta\sqrt t}
\le C_rQ^2L_f^{-2r}.
\]
The last inequality follows since exponential decay in $\sqrt t$ dominates every fixed power of $t$ and $\ln t$.
Substituting this bound into~\eqref{eq:finitef} proves~\eqref{eq:fixedf}.

If the mean far-field patterns coincide, estimate~\eqref{eq:fixedf} holds for every $0<\delta_f\le\delta_0$.
Letting $\delta_f\downarrow0$ gives $f_1=f_2$.

For the bandwidth-dependent estimate, apply~\eqref{eq:finitef} for $L\ge1$.
If $K\ge\ln|\ln\delta_f|$, taking the cutoff to be $K$ proves~\eqref{eq:incf}.
Otherwise, use~\eqref{eq:fixedf} and add the nonnegative data term to obtain~\eqref{eq:incf}.
\end{proof}

\subsection{Stability of the random source strength}
We first estimate the difference of the random source strengths from the covariance data over a finite frequency band.
\begin{proposition}\label{prop:local}
For $L\ge1$, we have
\begin{equation}\label{eq:local}
\norm{\mu_1-\mu_2}_{L^\infty}\le C(\varepsilon_L+ML^{-\gamma}).
\end{equation}
\end{proposition}
\begin{proof}
Set $a=c^2(\mu_1-\mu_2)$.
Choose $\gamma<\alpha<\min(1,s-3/2)$ and an even integer $n\ge2$ such that $n\alpha/(n+\alpha)\ge\gamma$.
By Sobolev embedding, $\norm a_{L^\infty}+|a|_{C^{0,\alpha}}\le CM$.
Here $|a|_{C^{0,\alpha}}=\sup_{x\ne y}|a(x)-a(y)|/|x-y|^\alpha$ is the H\"older seminorm on $\R^3$.
Let $\psi\in C_c^\infty(B_1)$ be nonzero and define the normalized localized test function
\[
v_{x,\rho}(y)=\frac{\rho^{-3/2}\psi((y-x)/\rho)}
{\norm{\rho^{-3/2}\psi((\cdot-x)/\rho)}_{\HH}},
\qquad x\in D,\quad0<\rho\le1.
\]
Since the coefficients of $-c^2\Delta$ and their derivatives are bounded and $n$ is even, differentiation and scaling give
\[
\norm{v_{x,\rho}}_{\HH}=1,\qquad
\norm{(-c^2\Delta)^{n/2}v_{x,\rho}}_{\HH}\le C_n\rho^{-n}.
\]
It follows from the spectral theorem that
\begin{equation}\label{eq:packet}
\norm{(I-\Pi_L)v_{x,\rho}}_{\HH}\le C_n(L\rho)^{-n}.
\end{equation}
Since $|v_{x,\rho}|^2c^{-2}dy$ has total mass one and support in $B_\rho(x)$, H\"older continuity gives
\[
|a(x)-\ip{M_av_{x,\rho}}{v_{x,\rho}}|\le CM\rho^\alpha.
\]
By~\eqref{eq:packet}, replacing $v_{x,\rho}$ by $\Pi_Lv_{x,\rho}$ in both arguments of the quadratic form changes its value by at most $C_nM(L\rho)^{-n}$.
Lemma~\ref{lem:compression} bounds the projected form in absolute value by $\varepsilon_L$.
Hence
\[
|a(x)|\le\varepsilon_L+CM\rho^\alpha+C_nM(L\rho)^{-n}.
\]
Choose $\rho=L^{-n/(n+\alpha)}$ to balance the two remainder terms.
Their sum is bounded by $CM L^{-n\alpha/(n+\alpha)}\le CM L^{-\gamma}$.
Since $|\mu_1-\mu_2|\le c_-^{-2}|a|$ and the continuous function $\mu_1-\mu_2$ vanishes outside $D$, estimate~\eqref{eq:local} follows.
\end{proof}

For fixed observation directions, define
\begin{equation}\label{eq:holcov}
\mathscr C_{\theta,\varphi}(z,w)=\int_D(\mu_1-\mu_2)(y)\Gamma(\theta,y;z)
\overline{\Gamma(\varphi,y;\overline w)}\,dy.
\end{equation}
The function is jointly holomorphic when both frequencies lie in the strip~\eqref{eq:strip}, equals $\mathcal C$ at real frequencies, and satisfies
\begin{equation}\label{eq:covgrowth}
|\mathscr C_{\theta,\varphi}(z,w)|\le CM(1+|z|)^3(1+|w|)^3.
\end{equation}
The reflected conjugate in the second factor is holomorphic in $w$.

\begin{proposition}\label{prop:covcont}
Suppose that~\eqref{eq:small} holds with $0<\delta_\mu<1$.
There exist $a,b,\beta,C>0$, independent of the source strengths and the data error, such that
\begin{equation}\label{eq:covcont}
\sup_{\theta,\varphi\in\sphere,\ 0<k,\ell\le L}|\mathcal C(\theta,\varphi;k,\ell)|
\le CM\exp\big(2bL-\beta|\ln\delta_\mu|e^{-2aL}\big),\qquad L\ge1.
\end{equation}
\end{proposition}
\begin{proof}
By~\eqref{eq:covgrowth}, we may fix $b>0$ and take $B\ge1$ sufficiently large so that
\[
H(z,w)=\frac{e^{-b(z+w)}}{BM}\mathscr C_{\theta,\varphi}(z,w)
\]
is bounded by one when both variables lie in the strip and by $\delta_\mu$ on $I^2$.
For fixed $\ell\in I$, Lemma~\ref{lem:propagation} in $z$ gives, for every real $k\ge0$, the bound $|H(k,\ell)|\le\delta_\mu^{\omega(k)}$.
Fixing $k$ and applying the lemma in $w$ with this bound on $I$, we obtain
\[
|H(k,\ell)|\le\delta_\mu^{\omega(k)\omega(\ell)}
\le\exp\big(-c_0^2|\ln\delta_\mu|e^{-a_0(k+\ell)}\big).
\]
The estimate follows by multiplying by $BM e^{b(k+\ell)}$ and using $k+\ell\le2L$.
\end{proof}

\begin{proof}[Proof of~\eqref{eq:fixedmu} and~\eqref{eq:incmu}]
Combining Lemma~\ref{lem:compression}, Proposition~\ref{prop:local}, and Proposition~\ref{prop:covcont}, we obtain, for every $L\ge1$,
\begin{equation}\label{eq:balance}
\norm{\mu_1-\mu_2}_{L^\infty}\le CM\left(
L^3\exp(2bL-\beta|\ln\delta_\mu|e^{-2aL})+L^{-\gamma}\right).
\end{equation}
Set $t=|\ln\delta_\mu|$ and choose
\begin{equation}\label{eq:choicemu}
L_\mu=\frac1{4a}\ln t.
\end{equation}
For sufficiently small $\delta_0$, the cutoff is at least one, and the first term in~\eqref{eq:balance} satisfies
\[
L_\mu^3t^{b/(2a)}e^{-\beta\sqrt t}\le C L_\mu^{-\gamma}.
\]
Substitution into the preceding bound proves~\eqref{eq:fixedmu}.

If the covariance data coincide, estimate~\eqref{eq:fixedmu} holds for every $0<\delta_\mu\le\delta_0$.
Letting $\delta_\mu\downarrow0$ yields $\mu_1=\mu_2$.

For the bandwidth-dependent estimate, apply~\eqref{eq:local} for $L\ge1$.
If $K\ge\ln|\ln\delta_\mu|$, take the cutoff to be $K$.
Otherwise, use~\eqref{eq:fixedmu} and add the nonnegative data term.
This proves~\eqref{eq:incmu}.
\end{proof}

\section{Explicit reconstruction and sampling}\label{sec:reconstruction}
In this section, we propose a reconstruction method for the stochastic inverse source problem. To do so,
we derive weighted backprojection formulas for the source mean $f$ and strength $\mu$ and estimate the reconstruction errors.
For each reference direction $\eta$ and fixed reference frequency $k_*$, the covariance $C(\theta,\eta;k,k_*)$, as a function of $(\theta,k)$, is the far-field pattern generated by the deterministic source $\mu(x)\overline{\Gamma(\eta,x;k_*)}$.
Applying the spectral identity and integrating over $\eta$, we obtain the weight for reconstructing the source strength.

\subsection{Backprojection from far-field data}
Fix a reference frequency $k_*>0$ and define
\begin{equation}\label{eq:referenceweight}
w_*(x)=c(x)^2\int_{\sphere}|\Gamma(\eta,x;k_*)|^2\,d\eta.
\end{equation}
For $\varphi\in L^2(D)$, define
\begin{equation}\label{eq:slicemap}
C_{\varphi,*}(\theta,\eta,k)=\int_D \varphi(x)\Gamma(\theta,x;k)
\overline{\Gamma(\eta,x;k_*)}\,dx.
\end{equation}
We write $C_*=C_{\mu,*}$ for the measured cross-frequency covariance with its second frequency fixed at $k_*$.

\begin{lemma}\label{lem:weightpositive}
The function $w_*$ is smooth and satisfies $0<w_-\le w_*(x)\le w_+$ on $D$ for some positive constants.
Moreover,
\begin{equation}\label{eq:sliceisometry}
\frac2\pi\int_0^\infty k^2\|C_{\varphi,*}(\cdot,\cdot,k)\|_{L^2(\sphere\times\sphere)}^2\,dk
=\int_Dw_*(x)|\varphi(x)|^2\,dx.
\end{equation}
\end{lemma}
\begin{proof}
If $w_*(x_0)=0$, the point-source field $G(\cdot,x_0;k_*)$ has a vanishing far-field pattern.
Rellich's lemma and unique continuation imply that the field vanishes away from $x_0$, contradicting the singularity of the Green function at $x_0$.
Mixed reciprocity and interior elliptic regularity give the smoothness.
The bounds follow from the compactness of $D$.
For each $\eta$, apply~\eqref{eq:flux} with $g=c^2\varphi\overline{\Gamma(\eta,\cdot;k_*)}$, so that $u_g^\infty(\theta,k)=C_{\varphi,*}(\theta,\eta,k)$.
Integrating over $\eta$ gives~\eqref{eq:sliceisometry}.
\end{proof}

\begin{remark}\label{rem:finitereference}
Let $\eta_a\in\sphere$ and $\rho_a>0$, $a=1,\ldots,J$, be the nodes and weights of a quadrature rule for the integral over the reference directions.
If the corresponding weight
\[
w_*(x)=c(x)^2\sum_{a=1}^J\rho_a|\Gamma(\eta_a,x;k_*)|^2
\]
has a positive lower bound on $D$, the isometry and reconstruction formulas hold with this rule used in both the data norm and the backprojection.
A single reference direction is therefore sufficient for reconstruction when its kernel is nonzero throughout $D$.
\end{remark}

For $0<K\le\infty$, let $\ZZ_K$ and $\YY_K$ be the weighted $L^2$ spaces on $(0,K)\times\sphere$ and $(0,K)\times\sphere\times\sphere$, respectively, with norms
\begin{equation}\label{eq:recondatanorms}
\|v\|_{\ZZ_K}^2=\frac2\pi\int_0^Kk^2\|v(\cdot,k)\|_{L^2(\sphere)}^2\,dk,
\quad
\|H\|_{\YY_K}^2=\frac2\pi\int_0^Kk^2\|H(\cdot,\cdot,k)\|_{L^2(\sphere\times\sphere)}^2\,dk.
\end{equation}
Define
\begin{equation}\label{eq:meanbackprojection}
(\meanback_Kv)(x)=\frac{2}{\pi c(x)^2}
\int_0^Kk^2\int_{\sphere}v(\theta,k)\overline{\Gamma(\theta,x;k)}\,d\theta\,dk
\end{equation}
and
\begin{equation}\label{eq:covbackprojection}
(\covback_KH)(x)=\frac{2}{\pi w_*(x)}\int_0^Kk^2
\int_{\sphere\times\sphere}H(\theta,\eta,k)
\overline{\Gamma(\theta,x;k)}\Gamma(\eta,x;k_*)\,d\theta\,d\eta\,dk.
\end{equation}
To reconstruct the source mean, we rescale the far-field data, apply the adjoint spectral transform, and divide by $c^2$.
For the source strength, the backprojection integrates the cross-frequency covariance over the reference directions and divides by $w_*$.
Integrals over all positive frequencies are understood as limits in $L^2(D)$.
We define the pointwise projection onto the nonnegative real axis by $\Pi_+v=\max(0,\re v)$.

\begin{theorem}\label{thm:explicit}
Let $f\in L^2(D)$ and let $\mu\in L^\infty(D)$ be nonnegative, with both functions extended by zero to the whole space.
Then
\begin{equation}\label{eq:fullinversion}
f=\lim_{K\to\infty}\meanback_K(\E U),\qquad
\mu=\lim_{K\to\infty}\covback_KC_*
\quad\hbox{in }L^2(D).
\end{equation}
Suppose in addition that, for $r>0$ and an integer $m\ge2$,
\[
\|(I-c^2\Delta)^{r/2}(c^2f)\|_{\HH}\le Q,\qquad
\|\mu\|_{H^m(\R^3)}\le M_m.
\]
If $\|e_f\|_{\ZZ_K}\le\eta_f$ and $\|e_\mu\|_{\YY_K}\le\eta_\mu$, the perturbed reconstructions
\begin{equation}\label{eq:truncatedestimators}
f_K=\meanback_K(\E U+e_f),\qquad
\mu_K=\Pi_+\covback_K(C_*+e_\mu)
\end{equation}
satisfy, for $K\ge1$,
\begin{equation}\label{eq:expliciterrors}
\|f_K-f\|_{L^2(D)}\le C(\eta_f+QK^{-r}),\qquad
\|\mu_K-\mu\|_{L^2(D)}\le C(\eta_\mu+M_mK^{-m}).
\end{equation}
The constant $C$ may depend on $c,D,k_*,r,m$, but is independent of $K$, the sources, and the perturbations.
\end{theorem}
\begin{proof}
It follows from the intertwining identities that
\begin{equation}\label{eq:meanprojector}
\meanback_K(\E U)=c^{-2}\Pi_K(c^2f).
\end{equation}
Since $\Pi_K\to I$ strongly on $\HH$, we obtain the inversion formula for every $f\in L^2(D)$.
The graph-norm bound gives $\|(I-\Pi_K)c^2f\|_{\HH}\le QK^{-r}$.
The isometry of $\FF$ also gives $\|\meanback_Kv\|_{L^2(D)}\le C\|v\|_{\ZZ_K}$.

For $H\in\YY_K$ and $\varphi\in L^2(D)$, Fubini's theorem gives
\[
\langle\covback_KH,\varphi\rangle_{L^2(D,w_*dx)}
=\langle H,C_{\varphi,*}\rangle_{\YY_K}.
\]
By the Cauchy--Schwarz inequality and~\eqref{eq:sliceisometry},
\[
\|\covback_KH\|_{L^2(D,w_*dx)}\le\|H\|_{\YY_K}.
\]
For $H\in\YY_\infty$ and $L>K$, this bound implies
\[
\|\covback_LH-\covback_KH\|_{L^2(D,w_*dx)}
\le\|\mathbf1_{\{K<k<L\}}H\|_{\YY_\infty}.
\]
The right-hand side tends to zero as $K\to\infty$, so we define
\[
\covback_\infty H=\lim_{K\to\infty}\covback_KH
\quad\hbox{in }L^2(D,w_*dx),\qquad
\|\covback_\infty H\|_{L^2(D,w_*dx)}\le\|H\|_{\YY_\infty}.
\]
Passing to the limit in the pairing identity yields
\[
\langle\covback_\infty H,\varphi\rangle_{L^2(D,w_*dx)}
=\langle H,C_{\varphi,*}\rangle_{\YY_\infty}.
\]
Polarizing~\eqref{eq:sliceisometry}, we obtain $\covback_\infty C_{\varphi,*}=\varphi$.
The lower bound on $w_*$ gives $\|\covback_KH\|_{L^2(D)}\le C\|H\|_{\YY_K}$, uniformly in $K$.
Taking $H=C_*=C_{\mu,*}$ proves the second inversion formula.

For the frequency truncation error, let $F(x,\eta)=\mu(x)\overline{\Gamma(\eta,x;k_*)}$ be an $L^2(\sphere)$-valued source.
The Hilbert tensor product $\HH\otimes L^2(\sphere)$ is identified with $L^2(\sphere;\HH)$, and the operator powers below act only on $x$.
The spectral theorem gives
\begin{equation}\label{eq:momentidentity}
\frac2\pi\int_0^\infty k^{2m+2}\|C_*(\cdot,\cdot,k)\|_{L^2(\sphere\times\sphere)}^2\,dk
=\|(-c^2\Delta)^{m/2}(c^2F)\|_{\HH\otimes L^2(\sphere)}^2\le CM_m^2.
\end{equation}
Since the acoustic operator is a smooth compactly supported perturbation of $-\Delta$, we have $D((I-c^2\Delta)^{m/2})=H^m(\R^3)$ with equivalent norms.
Multiplication by the smooth reference kernel is bounded on this space, uniformly in the reference direction.
This proves the last inequality in~\eqref{eq:momentidentity} and gives $\|\mathbf1_{\{k>K\}}C_*\|_{\YY_\infty}\le CM_mK^{-m}$.
Combining the truncation bounds with the boundedness of the backprojections gives~\eqref{eq:expliciterrors} before projection.
Since $\mu\ge0$, the pointwise inequality $|\Pi_+v-\mu|\le|v-\mu|$ gives the same error bound after projection.
\end{proof}

\subsection{Sample mean and covariance}
Let $U_n(\theta,k)$ be the far-field pattern in the $n$-th experiment, where $n=1,\ldots,N$ and $N\ge2$.
Different experiments are independent, and the measurements at all frequencies in one experiment share the same source realization.
We approximate the far-field statistics by the sample mean and cross-frequency sample covariance
\begin{equation}\label{eq:empiricalmean}
\overline U_N(\theta,k)=\frac1N\sum_{n=1}^NU_n(\theta,k)
\end{equation}
and
\begin{equation}\label{eq:empiricalcovariance}
\widehat C_{*,N}(\theta,\eta,k)=\frac1{N-1}\sum_{n=1}^N
(U_n(\theta,k)-\overline U_N(\theta,k))
\overline{(U_n(\eta,k_*)-\overline U_N(\eta,k_*))}.
\end{equation}
Using the observations paired by the sample index $n$, we reconstruct the source mean and strength by
\begin{equation}\label{eq:sampleestimators}
\widehat f_{K,N}=\meanback_K\overline U_N,\qquad
\widehat\mu_{K,N}=\Pi_+\covback_K\widehat C_{*,N}.
\end{equation}
For real $f$, we take the real part of $\widehat f_{K,N}$.

\begin{lemma}\label{lem:bandenergy}
For $K\ge1$,
\begin{equation}\label{eq:bandenergy}
\sup_{x\in D}\int_0^Kk^2\|\Gamma(\cdot,x;k)\|_{L^2(\sphere)}^2\,dk\le CK^3.
\end{equation}
\end{lemma}
\begin{proof}
For the Dirichlet form of $-c^2\Delta$ with measure $c^{-2}dx$, the three-dimensional Nash inequality gives the heat kernel estimate $p_t(x,x)\le Ct^{-3/2}$; see~\cite{CKKW}.
Here $p_t(x,y)$ is the integral kernel of $e^{tc^2\Delta}$ with respect to the measure $c(y)^{-2}\,dy$.
By~\eqref{eq:intertwine}, the diagonal of the integral kernel of $\Pi_K$ with respect to $c^{-2}dy$ is
\[
e_K(x,x)=\frac2\pi\int_0^K k^2\|\Gamma(\cdot,x;k)\|^2\,dk.
\]
Using $\mathbf1_{\{0\le\lambda\le K^2\}}\le e\,e^{-\lambda/K^2}$, we obtain $e_K(x,x)\le e\,p_{K^{-2}}(x,x)\le CK^3$.
\end{proof}

\begin{theorem}\label{thm:sampling}
Under the graph-norm and Sobolev bounds in Theorem~\ref{thm:explicit}, set $M_1=\|\mu\|_{L^1(D)}$.
Then
\begin{align}
\E\|\widehat f_{K,N}-f\|_{L^2(D)}^2
&\le C\left(Q^2K^{-2r}+\frac{M_1K^3}{N}\right),\label{eq:samplemeanrate}\\
\E\|\widehat\mu_{K,N}-\mu\|_{L^2(D)}^2
&\le C\left(M_m^2K^{-2m}+\frac{M_1^2K^3}{N-1}\right).
\label{eq:samplemurate}
\end{align}
The constants are independent of $K$ and $N$.
\end{theorem}
\begin{proof}
Let $Y=U-\E U$.
Using the independence of the experiments and the white-noise isometry, we obtain from~\eqref{eq:bandenergy}
\[
\E\|\overline U_N-\E U\|_{\ZZ_K}^2
=\frac2{\pi N}\int_D\mu(x)\int_0^Kk^2\|\Gamma(\cdot,x;k)\|^2\,dk\,dx
\le CM_1K^3/N.
\]
To estimate the sample covariance error, let $Z(k,\theta)=kY(\theta,k)$ and $V(\eta)=Y(\eta,k_*)$.
Set $\tau_Z=\E\|Z\|_{L^2((0,K)\times\sphere)}^2$ and $\tau_V=\E\|V\|_{L^2(\sphere)}^2$.
The corresponding total variances are bounded by $\tau_Z\le CM_1K^3$ and $\tau_V\le CM_1$.
Applying the Gaussian sampling identity in Appendix~\ref{app:sampling}, we obtain
\begin{equation}\label{eq:crosssampling}
\E\|\widehat C_{*,N}-C_*\|_{\YY_K}^2
\le\frac{C\tau_Z\tau_V}{N-1}\le\frac{CM_1^2K^3}{N-1}.
\end{equation}
The result follows by combining these estimates with Theorem~\ref{thm:explicit}.
\end{proof}

\begin{remark}
For $f\in L^2(D)$ and nonnegative $\mu\in L^\infty(D)$, we replace the power-law truncation terms by the omitted spectral energies.
These energies converge to zero by Theorem~\ref{thm:explicit}.
Thus, in the white-noise source model, $K_N\to\infty$ and $K_N^3/N\to0$ give mean-square consistency of both reconstructions.
\end{remark}

The estimates separate the frequency truncation error from the sampling error.
At fixed $N$, increasing $K$ reduces the truncation term and increases the sampling bound.
At fixed $K$, increasing $N$ reduces the sampling bound.
When centered Gaussian measurement errors are independent of the source and the observation and reference channels use independent errors, the sample cross-covariance remains unbiased.
The measurement errors contribute to the total variances $\tau_Z$ and $\tau_V$ and to the variance of the sample mean.
We include these measurement errors in the numerical experiments.

\begin{remark}\label{rem:2d}
The same energy identity and reconstruction formulas hold in two dimensions.
With the far-field convention $G(r\theta,y;k)=e^{ikr}r^{-1/2}\Gamma(\theta,y;k)+O(r^{-3/2})$, the free kernel and its far-field pattern are
\begin{equation}\label{eq:2dnormalization}
G_0(x,y;k)=-\frac i4H_0^{(1)}(k|x-y|),\qquad
\Gamma_0(\theta,y;k)=-\frac{e^{i\pi/4}}{\sqrt{8\pi k}}e^{-ik\theta\cdot y}.
\end{equation}
Here $H_0^{(1)}$ is the Hankel function of the first kind of order zero.
Since the boundary flux is $k\|u_g^\infty\|_{L^2(\mathbb S^1)}^2$, Stone's formula retains the weight $2k^2/\pi$.
Replacing $\sphere$ by $\mathbb S^1$ in~\eqref{eq:referenceweight}--\eqref{eq:covbackprojection} gives the two-dimensional backprojections.
The two-dimensional Nash inequality gives $CK^2$ in~\eqref{eq:bandenergy}, so the sampling terms in~\eqref{eq:samplemeanrate}--\eqref{eq:samplemurate} become $CM_1K^2/N$ and $CM_1^2K^2/(N-1)$, respectively.
\end{remark}

\section{Numerical experiments}\label{sec:numerics}
In this section, we describe the numerical implementation of the weighted backprojection formulas and present two-dimensional experiments for reconstructing the source mean and strength.
We study the effects of bandwidth and sample size using paired observations with 3\% relative root-mean-square (RMS) measurement noise.

\subsection{Implementation}
The known medium is
\begin{equation}\label{eq:numericalmedium}
c(x)^{-2}=1+0.35b(|x|/0.75),\qquad
b(t)=\begin{cases}\exp(1-(1-t^2)^{-1}),&0\le t<1,\\0,&t\ge1.\end{cases}
\end{equation}
The medium satisfies the Herglotz condition $\partial_r(r/c(r))>0$ and is therefore nontrapping.
We reconstruct the source functions in $B_{0.70}$.
Let $z_p$ be the quadrature nodes used to generate the data and $\omega_p$ the corresponding area weights.
At frequency $k_\ell$, the far-field data for the $n$-th experiment are generated by
\begin{equation}\label{eq:discretefield}
U_{n,j\ell}=\sum_p\omega_p\Gamma(\theta_{j\ell},z_p;k_\ell)f(z_p)
+\sum_p\sqrt{\omega_p\mu(z_p)}\Gamma(\theta_{j\ell},z_p;k_\ell)\xi_{np},
\end{equation}
where the $\xi_{np}$ are independent standard real Gaussian variables.
The same spatial source realization is used at all frequencies in each experiment.
We add circular complex Gaussian errors to the far-field measurements at the sampled directions according to
\begin{equation}\label{eq:measurementnoise}
\widetilde U_{n,j\ell}=U_{n,j\ell}
+0.03(\E|U_{n,j\ell}|^2)^{1/2}Z_{n,j\ell},\qquad
\E Z_{n,j\ell}=0,\quad\E|Z_{n,j\ell}|^2=1.
\end{equation}
Each $Z_{n,j\ell}$ has the representation $(X_{n,j\ell}+iY_{n,j\ell})/\sqrt2$, where $X_{n,j\ell}$ and $Y_{n,j\ell}$ are independent standard real Gaussian variables.
The measurement errors are independent of the source and independent across samples, measurement points, and channels.
The sample cross-covariance remains unbiased for the cross-covariance of the noiseless discrete data, and its variance includes the measurement noise contribution.
To compute the reconstructions, we evaluate~\eqref{eq:sampleestimators} by quadrature, take the real part of the reconstructed source mean, and project the reconstructed source strength onto nonnegative values.

For large $N$, we use Gaussian regression and the Bartlett factorization of a Wishart matrix~\cite[Chapter~3]{Muirhead1982} to draw the sample mean and sample cross-covariance of the discretized data from their exact joint finite-sample distribution.
Sharing the source and reference moments across frequencies preserves the joint dependence of the simulated statistics.
Here $N$ is the number of independent paired experiments represented by these statistics.

\subsection{Two-dimensional bandwidth and sampling effects}
We use the two-dimensional model with the normalization specified in Remark~\ref{rem:2d}.
The data generation and reconstruction employ distinct Cartesian midpoint grids, $80\times80$ and $81\times81$, in a square of side $1.5$.
The corresponding radial Lippmann--Schwinger solvers use 384 and 320 radial nodes, respectively.
At each frequency $k$, the observations are taken in $\max(48,4\lceil0.75k+8\rceil)$ uniformly spaced directions.
The reference frequency is $k_*=3$, and 48 reference directions are used.
For frequency integration, we use the composite Simpson rule with step size $0.5$ and cutoffs $K=8,16,28,40,56,72,80,96$.

The five strength profiles are a single Gaussian peak, two Gaussian peaks, an A shape, an L shape with a separate disk, and a five-petal flower with a central hole.
The corresponding source means include nonnegative and sign-changing profiles.

\begin{figure}[!t]
\centering
\includegraphics[width=\linewidth,height=0.78\textheight,keepaspectratio]{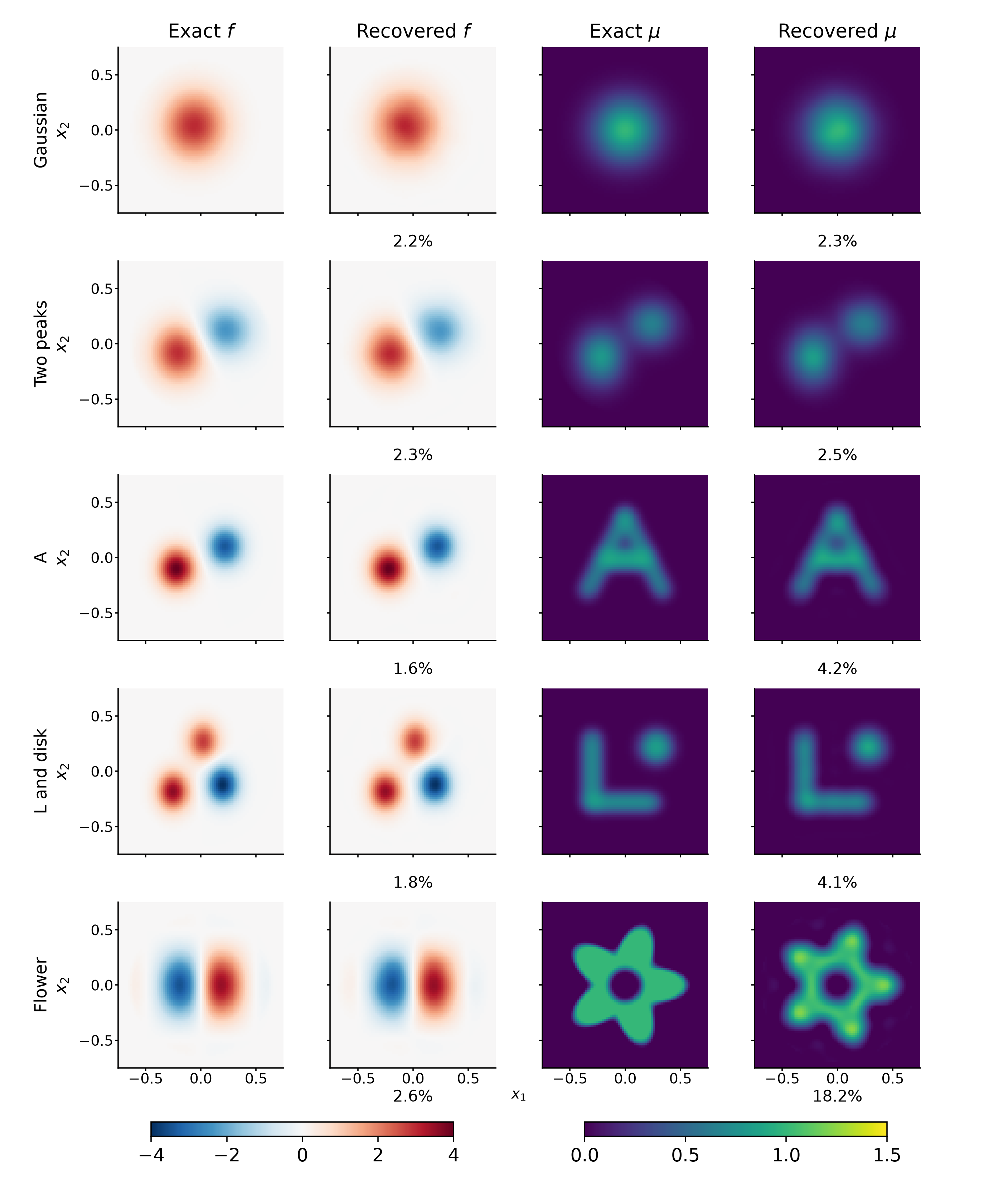}
\caption{Exact and reconstructed source means $f$ and strengths $\mu$ at $K=28$ and $N=32768$ with 3\% relative RMS measurement noise.
The percentages give the relative $L^2$ errors of the reconstructions from the first run.
All panels for $f$ share one color scale, and all panels for $\mu$ share another.}
\label{fig:reconstruction-k28}
\end{figure}

\begin{figure}[!t]
\centering
\includegraphics[width=\linewidth,height=0.78\textheight,keepaspectratio]{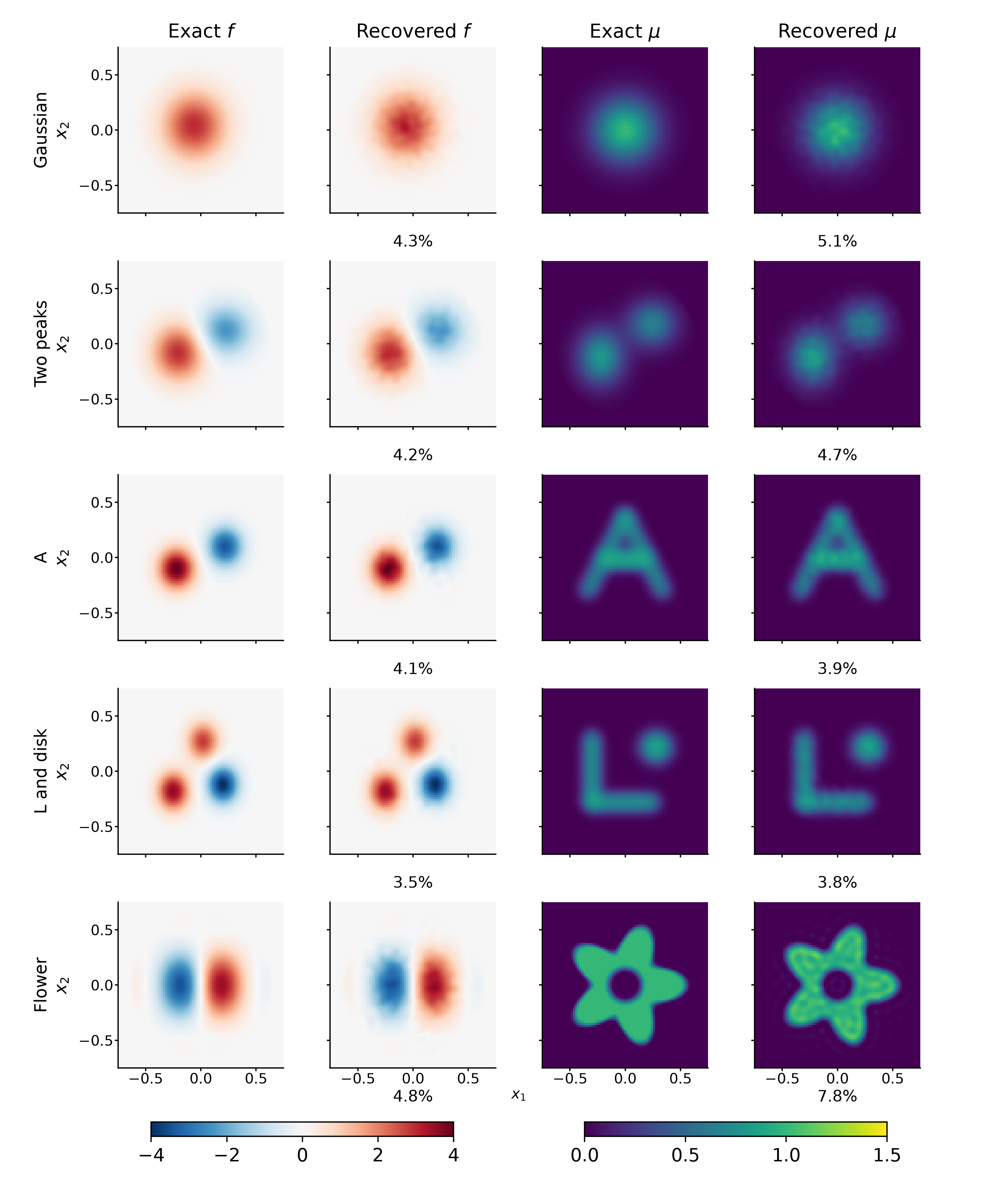}
\caption{Exact and reconstructed source means $f$ and strengths $\mu$ at $K=56$ and $N=32768$ with 3\% relative RMS measurement noise.
The percentages give the relative $L^2$ errors of the reconstructions from the first run.
All panels for $f$ share one color scale, and all panels for $\mu$ share another.}
\label{fig:reconstruction}
\end{figure}

Figures~\ref{fig:reconstruction-k28} and~\ref{fig:reconstruction} show the numerical reconstructions of the five source pairs at $K=28$ and $K=56$, respectively, with $N=32768$ and 3\% relative RMS measurement noise.
The relative $L^2(B_{0.70})$ errors are evaluated on the reconstruction grid.
The two figures use the same grids, panel layout, and color scales.
The flower strength has finer spatial features along the petal boundaries and around the central hole, which are less accurately resolved at $K=28$.
In the displayed run, the relative $L^2$ error of the reconstructed flower strength is $18.17\%$.
At $K=56$, the flower strength is reconstructed more accurately, although visible fluctuations remain.

\begin{figure}[!htbp]
\centering
\includegraphics[width=\linewidth]{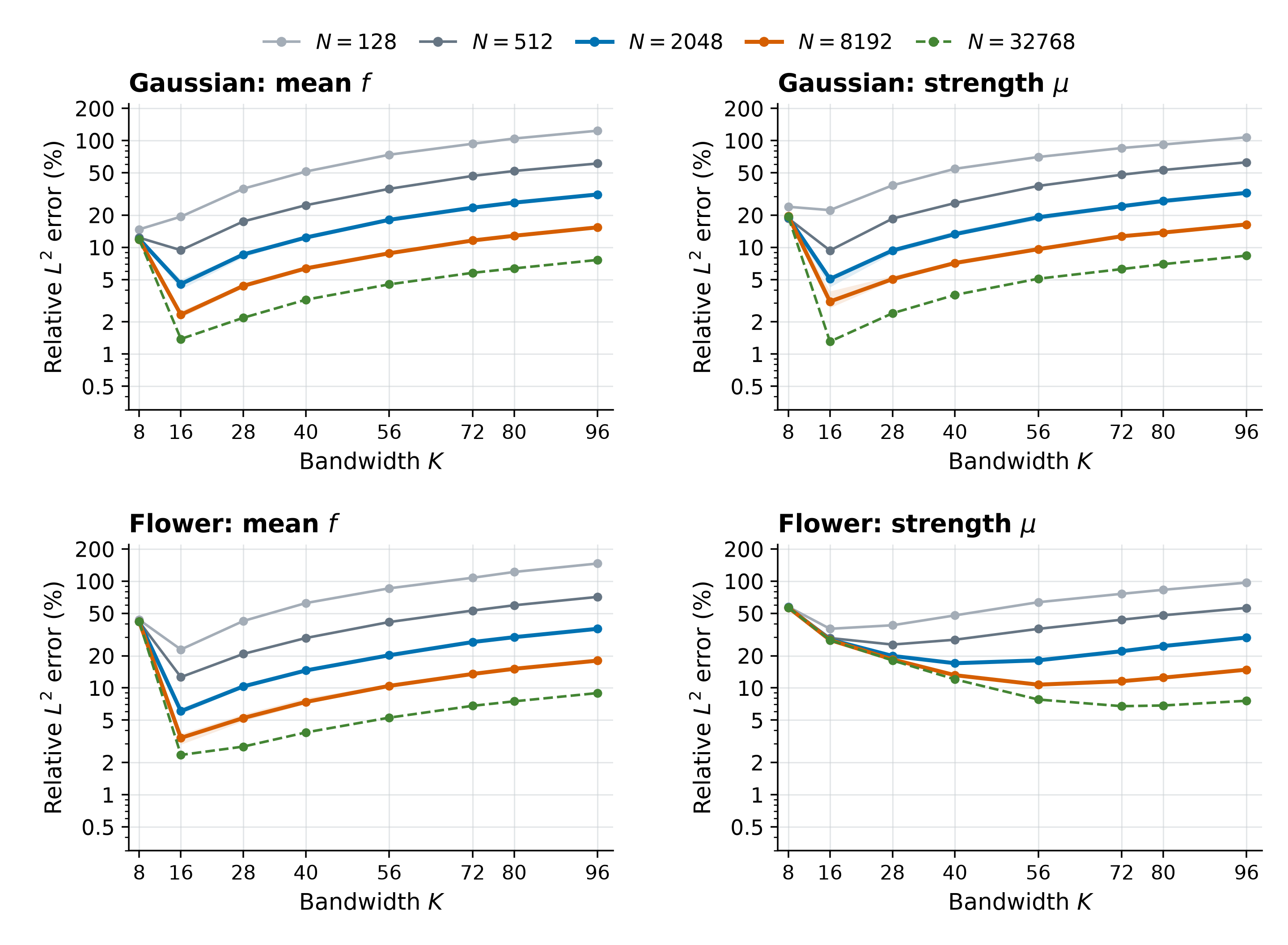}
\caption{Relative $L^2$ reconstruction errors at the eight cutoffs $K$ and five sample sizes $N$ for the two source pairs whose strength profiles are a single Gaussian peak and a five-petal flower with a central hole.
The curves show medians over ten independent runs with 3\% relative RMS measurement noise.
Shaded bands indicate the interquartile ranges for $N=2048$ and $N=8192$.
The vertical axes are logarithmic.}
\label{fig:medium-bandwidth}
\end{figure}

We examine the effects of bandwidth and sample size for the two source pairs whose strength profiles are a single Gaussian peak and a five-petal flower with a central hole.
For each pair, we consider $N=128,512,2048,8192,32768$ at all eight cutoffs.
We perform ten independent Monte Carlo runs for each source pair and sample size.
Within each run, the reconstructions at different cutoffs use the same sample mean and sample cross-covariance, restricted to the corresponding frequency bands.
Data sets for different values of $N$ are generated independently.
Figure~\ref{fig:medium-bandwidth} shows the resulting median relative $L^2$ reconstruction errors.

\clearpage

 Figure~\ref{fig:medium-bandwidth} shows that, at each tested cutoff $K\ge16$, the median reconstruction errors decrease as $N$ increases.
For $N\ge512$, the median errors initially decrease as $K$ increases and grow again at larger cutoffs.
For the flower strength, the decrease extends to larger cutoffs as more samples become available.
These results illustrate the balance between frequency truncation and sampling errors.

\section{Conclusion}\label{sec:conclusion}
We have studied an inverse random source problem in a known inhomogeneous medium using multifrequency far-field statistics.
For a known nontrapping wave speed and under a priori regularity bounds, we have established conditional $L^2$ stability for the source mean and $L^\infty$ stability for the source strength from their respective far-field statistics on a fixed frequency interval.
We also derive weighted backprojection formulas for the source mean from the mean far-field pattern and for the source strength from the cross-frequency covariance.
The mean-square error bounds separate the contributions of frequency truncation and sampling, which yields efficient reconstruction methods for the stochastic inverse source problem.
Numerical experiments demonstrate the efficiency of the proposed reconstruction methods and show the effects of bandwidth and sample size.

\appendix
\section{Gaussian sampling identities}\label{app:sampling}
Let $Z$ and $V$ be centered, jointly complex Gaussian random elements in two Hilbert spaces, with coordinates $Z_i$ and $V_j$ in their respective orthonormal bases.
Write $C_{ZV}$ for the cross-covariance operator from the Hilbert space of $V$ to that of $Z$, with matrix entries $(C_{ZV})_{ij}=\E(Z_i\overline{V_j})$.
For Gaussian variables generated by real white noise, the nonconjugated cross-covariances $\E(Z_iV_j)$ need not be zero.
Denote the total variances by $\tau_Z=\E\|Z\|^2$ and $\tau_V=\E\|V\|^2$.
Let $\widehat C_{ZV,N}$ be the unbiased sample cross-covariance operator computed from $N\ge2$ independent copies of $(Z,V)$.
Then
\begin{equation}\label{eq:exactcrossvariance}
\E\|\widehat C_{ZV,N}-C_{ZV}\|_{\rm HS}^2
=\frac{\tau_Z\tau_V+\sum_{i,j}|\E(Z_iV_j)|^2}{N-1}
\le\frac{2\tau_Z\tau_V}{N-1}.
\end{equation}
By Isserlis' fourth-moment formula~\cite{Isserlis1918}, applied to the real and imaginary parts, each pair of coordinates satisfies
\[
\E|Z_i\overline{V_j}-\E(Z_i\overline{V_j})|^2
=\E|Z_i|^2\E|V_j|^2+|\E(Z_iV_j)|^2.
\]
Summing this identity over the coordinates gives the numerator in~\eqref{eq:exactcrossvariance}.
A real orthogonal change of the sample coordinates separates the sample mean from $N-1$ independent Gaussian contrasts.
The unbiased covariance is the average of the corresponding $N-1$ cross-products, giving the denominator.
This is the normal-sample decomposition for the Wishart distribution; see~\cite[Chapter~3]{Muirhead1982}.
Finite-dimensional projections and monotone convergence give the identity for Hilbert-space variables.
The final inequality follows from the Cauchy--Schwarz bound $\sum_{i,j}|\E(Z_iV_j)|^2\le\tau_Z\tau_V$.

Taking a single centered Gaussian field $Y$, with $C_Y=C_{YY}$, $\widehat C_{Y,N}=\widehat C_{YY,N}$, and $\tau_Y=\E\|Y\|^2$, we similarly obtain
\begin{equation}\label{eq:exactcovariancevariance}
\E\|\widehat C_{Y,N}-C_Y\|_{\rm HS}^2
=\frac{\tau_Y^2+\sum_{i,j}|\E(Y_iY_j)|^2}{N-1}.
\end{equation}
For~\eqref{eq:crosssampling}, the Hilbert spaces are $L^2((0,K)\times\sphere)$ and $L^2(\sphere)$, and the Gaussian variables are $Z(k,\theta)=kY(\theta,k)$ and $V(\eta)=Y(\eta,k_*)$.
The factor $2/\pi$ in the data norm is included in the constant.

\section*{Acknowledgements}
X. Xu is partially supported by the National Key Research and Development Program of China (No. 2024YFA1012303), National Natural Science Foundation of China (No. 12525112), and the Open Research Project of Innovation Center of Yangtze River Delta, Zhejiang University.
Y. Zhao is supported by National Natural Science Foundation of China (No. 12371421).

\section*{Conflict of Interest Statement}
The authors declare that they have no conflict of interest.

\section*{Data Availability Statement}
The numerical experiments use synthetic data generated using the models and parameters described in this article.
The data supporting the numerical results are available from the corresponding author upon reasonable request.

\end{document}